\documentclass[a4paper, 11pt, twoside, reqno]{amsart} \usepackage[a4paper,inner=2.5cm,outer=2.5cm,top=3cm,bottom=3cm]{geometry}
\usepackage{amsmath,amscd}
\usepackage{amssymb}

\usepackage{amsthm}
\usepackage{comment}
\usepackage{mathrsfs}
\usepackage{graphicx, xcolor}
\usepackage{xcolor}
\usepackage{mathtools}
\usepackage[normalem]{ulem}
\usepackage[ocgcolorlinks, linkcolor=blue,citecolor=red,urlcolor=cyan]{hyperref}
\usepackage{bm}
\usepackage{bbm}
\usepackage{url}
\usepackage[utf8]{inputenc}
\usepackage{mathtools,amssymb}
    \usepackage{amssymb}
\usepackage{tikz}
\usepackage{dsfont}
\usepackage{relsize}
\usepackage{url}
\usepackage{xcolor}
\usepackage{graphicx}
\usepackage{mathrsfs}
\usepackage[shortlabels]{enumitem}
\usepackage{lineno}
\usepackage{amsmath}
\usepackage{enumitem}
\usepackage{amsthm} 
\usepackage{verbatim}
\usepackage{dsfont}
\usepackage[utf8]{inputenc}
\usepackage{tikz}
\numberwithin{equation}{section}

\allowdisplaybreaks
 
 \mathtoolsset{showonlyrefs}
\theoremstyle{definition}
\newtheorem{step}{Step}
\newtheorem{theorem}{Theorem}[section]

\newtheorem{example}[theorem]{Example}

\allowdisplaybreaks[0]
\usepackage{blindtext}
\usepackage{cleveref}
\renewcommand{\O}{\Omega}
\renewcommand{\o}{\omega}

\newcommand{\Rb}{\mathbb{R}}

\usepackage{cellspace}
\usepackage{makecell}
\title[Reconstruction of coefficient in dynamical Schr{\"o}dinger equation]{Reconstruction of time-dependent coefficients in a semilinear dynamical Schr{\"o}dinger equation}
\author[Kumar, Nakamura and Vashisth]{Parveen Kumar$^{\diamond}$, Gen Nakamura$^{*}$, and Manmohan Vashisth$^{\diamond}$}
\address{{$^{\diamond}$ Department of Mathematics, Indian Institute of Technology Ropar, Rupnagar,  Punjab - 140001, India.
		\newline
		\indent E-mail:{\tt\  pswamifca@gmail.com, manmohanvashisth@iitrpr.ac.in}}}
        \address{$^{*}$Department of Mathematics, Hokkaido University, Sapporo 060-0810, Japan.
	\newline\indent E-mail:{\tt \ nakamuragenn@gmail.com}}

\DeclareMathOperator{\supp}{supp} 
\begin{document}
	\begin{abstract}
In the present manuscript, we study an inverse problem related to a semilinear dynamical Schr{\"o}dinger  equation with lower order terms,   in a bounded domain of $\Rb^{1+n},n\geq 2$. Our focus is on determination of the time-dependent coefficients appearing in the aforementioned equation, from the boundary measurements of the solutions. More precisely,  we establish  the {pointwise reconstruction} formulae for determining the time-dependent coefficients of linear and nonlinear terms from the knowledge of Dirichlet-to-Neumann map. Since the concerned non-linear Schr\"odinger equation
possesses a trivial solution, we linearize the equation around the trivial solution and    use the  asymptotic solutions (\textit{with concentrated amplitudes}) of the linearized problem for reconstructing the aforementioned coefficients.  To be more specific, we use first-order linearization to reconstruct vector and scalar  potentials associated with the coefficients of linear terms and the higher-order linearization technique is used  to reconstruct coefficients of nonlinearity. The nonlinear equation considered in this manuscript can be seen as a generalization of  the Gross-Pitaevskii equation (GPE), which is employed to describe the dynamics of dilute Bose-Einstein condensates (BEC). 

\vspace{.2cm}
		

	\end{abstract}
	\maketitle
   \section{Introduction}
 \subsection{Mathematical formulation and problem of interest}
The present manuscript is concerned with  an inverse problem associated with an initial boundary value problem (IBVP)  for semilinear dynamical Schr{\"o}dinger equation. 
More precisely, for $m\geq 2$ an integer,  we consider the following IBVP associated with the aforementioned equation by 
\begin{align}\label{govern_eqn}
\begin{aligned}
\begin{cases}
\left(\mathrm{i}\partial_t  +\sum_{j=1}^{n}\left(\partial_j+\mathrm{i}A_j(t,x)\right)^2   + q(t,x)\right)u(t,x)=B(t,x,u,\overline{u}),~(t,x)\in Q, \\
u(0,x) = 0,\quad x\in\Omega,\\
u(t,x) = f(t,x), \quad (t,x)\in \Sigma:=(0,T)\times\partial\Omega,
\end{cases}
\end{aligned}
\end{align}
where $T>0$ be fixed, $Q:=(0,T)\times\Omega$ with $\Omega\subset \Rb^n~(n\geq 2)$ is a bounded domain with smooth boundary $\partial \Omega$ and magnetic potential    
$A(t,x):=(A_1(t,x),A_2(t,x),\cdots, A_n(t,x))$ be a real-vector field. Here we  denoted by  $\partial_t$ and $\partial_j$ ($1\leq j\leq n)$, the partial derivatives  with respect to time-variable $t\in (0,T)$ and  $x_j$, the $j$th component of spacial variable $x:=(x_1,x_2,\cdots,x_n)\in\Omega$, respectively. 
The non linear term $B(t,x,u,\overline{u})$ is given by 
\begin{align}\label{non linear term in govern equation} B(t,x,u,\overline{u}):=\sum_{a=1}^{m}r_a(t,x)u^a(t,x)(\overline{u})^{m-a}(t,x).
\end{align} 
Throughout  this manuscript, we assume that $A \in C_c^{\infty}(Q,\mathbb{R}^n)$ and $ q, r_a \in C_c^\infty(Q,\mathbb{R})$  for each integer $1\leq a\leq m$. Also, we denote the linear dynamical Schr\"odinger operator associated to \eqref{govern_eqn} by $\mathcal{P}_{A,q}$ and is given by 
\begin{align}\label{linear operator P}
\begin{aligned}
    \mathcal{P}_{A,q}&:= \mathrm{i}\partial_t  +\sum_{j=1}^{n}\left(\partial_j+\mathrm{i}A_j(t,x)\right)^2   + q(t,x)\\
    &=i\partial_t +\Delta+2iA\cdot \nabla+i(\nabla\cdot A)-\lvert A\rvert^2+q.
    \end{aligned} 
\end{align}
To discuss the well-posedness of the IBVP \eqref{govern_eqn}, we first introduce some function spaces. 
For $N=\O$, or $\partial\Omega$, and $s,t\geq 0$,   we define the Hilbert space $H^{s,t}((0,T)\times N)$ by 
\[H^{s,t}((0,T)\times N):= H^s(0,T;L^2(N))\cap L^2(0,T;H^t(N))\]
with the norm 
\[\lVert u \rVert_{H^{s,t}((0,T)\times N)} :=\lVert u \rVert_{H^s(0,T;L^2(N))}+\lVert u \rVert_{L^2(0,T;H^t(N))}.\]
Now if $s=t$, then from  (Proposition 2.3, \cite{Lions_vol2}), we can identify $H^{s,s}(Q)$ and $H^{s,s}(\Sigma)$ by $H^{s}(Q)$ and $H^{s}(\Sigma)$, respectively and their associated norms are equivalent.
The well-posedness of the forward problem \eqref{govern_eqn} follows from (Proposition $2.3$, \cite{lai2024partial}) for the small Dirichlet data $f \in \mathcal{E}_{\delta}(\Sigma)$ where $\mathcal{E}_{\delta}(\Sigma)$ is given by 
\begin{align}\label{definition of d delta}
\begin{aligned}\mathcal{E}_{\delta}(\Sigma):=\Big\{ f \in H^{2 \kappa+\frac{3}{2}}(\Sigma): \partial_t^\ell f(0,\cdot)=0  \mbox{ on } \partial \Omega, \mbox{ for  } \ell \leq  2 \kappa +1,~l\in \mathbb{N}  \mbox{ and } \lVert f \rVert_{H^{2 \kappa+(3/2)}(\Sigma)} \leq \delta\Big\}
\end{aligned}\end{align}
for $\delta>0$ sufficiently small. That is, there exists a unique solution
\[
u \in \mathcal{H}_0^{2 \kappa}(Q):=\{h\in H^{2 \kappa}(Q): \partial_t^\ell f(0,\cdot)=0  \mbox{ in } \Omega \mbox{ for  } \ell < 2 \kappa -1,l\in \mathbb{N} \}
\] to the IBVP \eqref{govern_eqn} which satisfies the following estimate
\[  \lVert u \rVert_{H^{2\kappa}(Q)}\leq C \lVert f \rVert_{H^{2\kappa+(3/2)}(\Sigma)},
\]
where $ \kappa> \frac{n+1}{2}$ is an integer. Also, we denote 
\begin{align}
    \widetilde{\mathcal{H}}^{2\kappa+\frac{3}{2}}_0(\Sigma):=\{g \in H^{2 \kappa+\frac{3}{2}}(\Sigma): \partial_t^\ell g(0,\cdot)=0  \mbox{ on } \partial \Omega, \mbox{ for  } \ell \leq 2 \kappa +1,~l\in \mathbb{N}
    \}.
\end{align}
  Based on the aforementioned well-posedness of \eqref{govern_eqn}, we define the Dirichlet-to-Neumann (DN) map  $\Lambda_{A,q,B}:\mathcal{E}_{\delta}(\Sigma) \rightarrow H^{2\kappa -\frac{3}{2}}(\Sigma)$, associated to \eqref{govern_eqn} by 
\begin{align}\label{DN map for gover eqn}
\begin{aligned}
\Lambda_{A,q,B}(f):=\partial_{\nu} u_{f} \Big|_{ \Sigma},\ \  f\in \mathcal{E}_{\delta}(\Sigma)
\end{aligned}
\end{align}
where $u_f$ is a solution to IBVP \eqref{govern_eqn} and $\partial_\nu:=\nu\cdot\nabla_x$, normal derivative with respect to the outward  unit  normal $\nu$ on $\partial\Omega$. 
In the current manuscript, we are interested in the following inverse problem associated with the IBVP \eqref{govern_eqn}.\\

{\bf{Problem of interest:}} To recover the coefficients $A$, $q$, and $r_a$ appearing in \eqref{govern_eqn} from the knowledge of the aforementioned DN map $\Lambda_{A,q,B}(f)$, measured for all $f\in \mathcal{E}_{\delta}(\Sigma)$.  More precisely, we prove Theorem \ref{Main theorem}, stated below as a  main result of this manuscript.
\subsection{Physical significance} The Schr{\"o}dinger equation, named after Erwin Schr{\"o}dinger is a building block of Quantum mechanics. It has a variety of applications, including quantum tunnelling, wave propagation in fibre optics, and more. The semilinear equation \eqref{govern_eqn} considered in  this manuscript has numerous applications in science and engineering, for example it is used to study physical processes such as plasma physics (like Langmuir wave dynamics \cite{khater2024langmuir}), nonlinear optics (like, optical pulse propagation \cite{akhmediev1994dynamics, akram2024study}), water waves (like, gravity waves on water layer \cite{hasimoto1972nonlinear, li2019nonlinear}), ocean surface waves (like, rogue wave/Killer waves \cite{ohta2012general}) and others. 
Also, we refer \cite{swiecicki2016schrodinger} for the connection between the semilinear equation \eqref{govern_eqn} and mean-field theories. 
Moreover, if we we choose  $B = r_2 u^2 \overline{u}$ and \(A = \mathbf{0}\), in  \eqref{govern_eqn}, then it reduces to  a well-known \textquotedblleft the Gross–Pitaevskii equation (GPE)\textquotedblright (see for example \cite{lassas2025coefficient} and references therein), which is employed to describe the dynamics of dilute Bose–Einstein condensates (BEC). Thus, the semilinear equation \eqref{govern_eqn} can be seen as a generalization of GPE. 
\subsection{Related works}   
The inverse problems related to  identification of  coefficients appearing in  Schr{\"o}dinger equations have been a  key object of study due to its appearance in various physical phenomena (see, for example  \cite{hasimoto1972nonlinear, akhmediev1994dynamics, khater2024langmuir}).  For instance, Sun in \cite{sun1993inverse} considered an inverse problem associated to unique determination of the time-independent magnetic field and electrical potential in a linear Schr{\"o}dinger equation from the DN map.  Bellassoud and Ferreira
 \cite{bellassoued2010stable} investigated a dynamical anisotropic  Schr{\"o}dinger equation with $A=0$ and derived 
 a H{\"o}lder-type stability for the determination of the potential in a manifold setup. In a periodic quantum waveguide, Choulli et al \cite{choulli2015stable} considered a linear  Schr{\"o}dinger equation with $A=0$ and proved a stable estimate for the scalar potential from the boundary measurement. 
For the dynamical  Schr{\"o}dinger equation, A{\"\i}cha  \cite{ben2017stability} established a stability estimate for the recovery of the magnetic field and electrical potential from the DN map.
Kian and Soccorsi  \cite{kian2019holder}, established  the H\"older-type stability in determining the electromagnetic potential $(A,q)$ in the  Schr{\"o}dinger equation and also showed that the unique recovery is achievable only if $\nabla\cdot A=0$ and  due to a natural gauge invariance, recovery of divergence free magnetic potential is optimal.
   For inverse problems related to uniqueness and stability in context of linear dynamical Schr{\"o}dinger equation, we refer to \cite{eskin2008inverse, avdonin2005dynamical, bellassoued2017stable, bellassoued2018inverse, kian2019holder} and references therein. In a manifold setup, for more works related to the determination of electrical and magnetic potentials, we refer to \cite{baudouin2007uniqueness, bellassoued2019stable, kian2019holder, krupchyk2023inverse}.

So far, we have mentioned the prior works related to recovery of coefficients appearing in linear static and dynamical Schr\"odinger equations. Next we mention, the existing results related to coefficients identification inverse problems for nonlinear PDEs.  
Motivated by the pioneer work of Isakov \cite{isakov1993uniqueness}, which introduced first time a linearization technique to deal with the inverse problem for nonlinear PDEs,  we are concerned with coefficient determination inverse problem related to a  semilinear dynamical  Schr{\"o}dinger equation. In this approach, we linearize the DN map and use the existing results for linear PDEs to solve the inverse problems related to  nonlinear PDEs. In \cite{lassas2021inverse}, the inverse problem for a semilinear elliptic equation with power-type nonlinearities is analyzed using a higher-order linearization method. In the context of nonlinear  Schr{\"o}dinger equation, Lai et al. \cite{lai2023partial} investigated a polynomial type nonlinearity in the absence of a magnetic potential.
By constructing geometric optics solutions that are based on Gaussian beam quasimodes, they achieved a global uniqueness from partial boundary measurements. Moreover, they establish a logarithmic stability estimate for the recovery of the nonlinear coefficient using the unique continuation principle. For the corresponding problem in a Riemannian manifold setting, Lasses et al. \cite{lassas2025coefficient} established the unique recovery of the coefficients from the source-to-solution map.
Later, for a more general nonlinear Schr{\"o}dinger equation, Lai et al. in \cite{lai2024partial} studied an inverse problem related to unique coefficient identification from the partial DN map.  In the context of reconstruction of the coefficients, Carstea et al. \cite{carstea2019reconstruction} considered a quasilinear elliptic operator of the divergence form and derived the reconstruction formulae for determining the  coefficients from the DN map. 
Recently, Bhardwaj et al.  \cite{bhardwaj2026reconstruction} established the reconstruction formulae for recovering  the potential and damping coefficients in a semilinear wave equation with a power-type nonlinearity. In addition, for the nonlinear fourth-order Schr{\"o}dinger equation, we refer the reader to \cite{MPS26}, where the uniqueness of the associated coefficients is established in both Euclidean and Riemannian manifold settings.
For more works on inverse problems related to   nonlinear PDEs,  which are closely related to the inverse problem considered in this manuscript,  we refer the reader   to \cite{lassas2021inverse, kang2002identification, krupchyk2020remark, krupchyk2021partial, carstea2019reconstruction, feizmohammadi2020inverse, Lai_Zhou, lassas2020partial}, to  \cite{choulli2018logarithmic, kaltenbacher2021simultaneous, kian2020recovery} and to \cite{nakamura2008inverse,chen2021detection, lassas2018inverse, kurylev2018inverse, nakamura2021inverse, bhardwaj2026reconstruction} for recovering the coefficients appearing in nonlinear elliptic,  parabolic and hyperbolic PDEs, respectively. 
\subsection{Statement of the main result} 
As mentioned earlier, we observe that from an inverse problem point of view, substantial progress has been made toward unique identification of the coefficients and stability aspects for both linear and nonlinear PDEs; however, few results are available for reconstructing the coefficients from boundary measurement data. Our main result stated below, provides the point-wise reconstruction formulae for determining the time-dependent coefficients of linear and nonlinear terms, appearing in semilinear dynamical Schr\"odinger equation \eqref{govern_eqn} from the knowledge DN map measured on $\Sigma$.  
\begin{theorem}\label{Main theorem}
 {\it    Let $\O$ be an open, bounded, and simply connected subset of $\mathbb{R}^n$ with smooth boundary $\partial \O$. Also, suppose that $A \in C_c^{\infty}(Q,\mathbb{R}^n)$ and  $ q,~r_a \in C_c^\infty(Q,\mathbb{R})$ for each integer $1\leq a\leq m$.  
Then, from the knowledge of DN map $\Lambda_{A,q,B}(f)$ where $f \in \mathcal{E}_{\delta}(\Sigma)$, we can uniquely reconstruct the coefficients $A$, $q$ and $r_a$ in $Q$ pointwise provided $\nabla\cdot A=0 \text{ in } Q$.}
\end{theorem}
Recently, in \cite{lai2024partial},  authors established a partial data unique determination result related to recovery of the  coefficients appearing in  nonlinear magnetic Schr{\"o}dinger equation, in which they proved the uniqueness of inverse problems when the coefficients are known in a small enough open neighborhood of the boundary of the spatial domain. Inspired by \cite{lai2024partial}, we  consider  the inverse problem related to finding  reconstruction formulae for determining the time-dependent coefficients appearing in \eqref{govern_eqn} from the boundary measurements of solution. Our result can  be viewed as a continuation of the analysis undertaken in \cite{lai2024partial}. To the best of our knowledge, our  work is the first result in the direction of reconstruction of time-dependent coefficients in semilinear dynamical Schr\"odinger  using the knowledge of DN map $\Lambda_{A,q,B}$.
\\

The rest of the manuscript is organized as follows.   \Cref{section 1} is devoted to proving the main result of this manuscript, where in Subsection \ref{sub1} we provide the reconstruction formulae for the vector potential $A$ and $q$ while in Subsection \ref{sub2} we provide the reconstruction formulas for the nonlinear coefficients. 
    \section{Proof of Theorem \ref{Main theorem}}\label{section 1}
    This section is devoted to proving the main result of this manuscript. First, we provide the reconstruction formulae for the vector and scalar potentials $A$ and $q$, respectively. Inspired by  prior works (see for example \cite{isakov1993uniqueness, nakamura2008inverse, bhardwaj2026reconstruction, lassas2021inverse} and references therein) related to  coefficients identification inverse  problems  for nonlinear PDEs, we use the first-order linearization technique  to recover the linear coefficients $A$ and $q$ and the reconstruction of coefficients  $r_i$ ($1\leq i\leq m$) of nonlinear terms, are established  via the higher-order linearization techniques along with asymptotic solutions to the linearized equation. We divide the present section into two subsections, in which first one consists of the proof related to reconstruction of the coefficients of linear terms and the second subsection consists of reconstruction of the coefficients of nonlinear terms. 
    \subsection{Reconstruction of coefficients associated with linear terms}\label{sub1}
In this subsection, we derive reconstruction formulae for finding vector and scalar potentials $A$ and  $q$ respectively. As mentioned above, following \cite{isakov1993uniqueness, nakamura2008inverse, bhardwaj2026reconstruction, lassas2021inverse}, we use  the first-order linearization of solution to the nonlinear problem  to recover these coefficients. Given $f_1, f_2,\cdots,f_m \in \widetilde{\mathcal{H}}^{2\kappa+\frac{3}{2}}_0(\Sigma)$,  choose   $\epsilon:=(\epsilon_1,\epsilon_2,\cdots,\epsilon_m)$ with $\epsilon_i\geq 0$ ($1\leq i\leq m$), such that $\epsilon f:=\sum_{i=1}^m \epsilon_if_i \in \mathcal{E}_{\delta}(\Sigma)$ and denote 
$u:=u_{\epsilon f}$, a solution to the following IBVP
\begin{align}\label{govern pde with boundary data ef}
\begin{cases}
\mathcal{P}_{A,q}u(t,x)=B(t,x,u,\overline{u}),~(t,x)\in Q, \\
u(0,x) = 0,\quad x\in\Omega,\\
u(t,x) = \epsilon f(t,x), \quad (t,x)\in \Sigma.
\end{cases}
\end{align}
Also, the Dirichlet-to-Neumann map is given by
\begin{align}\label{DN map for gover eqn 1}
\begin{aligned}
\Lambda_{A,q,B}(\epsilon f):=\partial_{\nu} u |_{ \Sigma}=\text{Known}, \ \mbox{for any $\epsilon f\in \mathcal{E}_{\delta}(\Sigma)$}
\end{aligned}
\end{align}
is well-defined, whenever 
 $u$ is a solution to \eqref{govern pde with boundary data ef}.

\subsubsection{First-order linearization} The first order linearization is done via differentiating \eqref{govern pde with boundary data ef} with respect to  $\epsilon_i$ ($1\leq i\leq m$) and substituting $\epsilon=0$.  Now after  differentiating \eqref{govern pde with boundary data ef} with respect to $\epsilon_i$ for $1\leq i \leq m$, we get
\begin{align*}
    \begin{cases}
        \mathrm{i}\partial_t (\partial_{\epsilon_i}u) +\sum_{j=1}^{n}\left(\partial_j+\mathrm{i}A_j\right)^2 \partial_{\epsilon_i}u + q(t,x)\partial_{\epsilon_i}u=\partial_{\epsilon_i}B(t,x,u,\overline{u}),~(t,x)\in Q,\\
        \partial_{\epsilon_i}u(0,x) = 0,\quad x\in\Omega,\\\partial_{\epsilon_i}u(t,x) = f_i(t,x), \quad (t,x)\in \Sigma,
    \end{cases}
\end{align*}
where \[\partial_{\epsilon_i}B(t,x,u,\overline{u})=\sum_{a=1}^mr_a\left( a\partial_{\epsilon_i}u^{a-1}(\overline{u})^{m-a}+(m-a)u^a (\overline{u})^{m-a-1}\right).
\]
Now due the well-posedness of the forward problem, we observe that $\displaystyle \partial_{\epsilon_i}B(t,x,u,\overline{u})|_{\epsilon=0}=0$. Hence, after evaluating the above equations at $\epsilon=0$, we obtain
\begin{align*}
    \begin{cases}
        \mathrm{i}\partial_t (\partial_{\epsilon_i}u\lvert_{\epsilon=0}) +\sum_{j=1}^{n}\left(\partial_j+\mathrm{i}A_j\right)^2 \partial_{\epsilon_i}u\lvert_{\epsilon=0}  + q(t,x)\partial_{\epsilon_i}u\lvert_{\epsilon=0}=0,~(t,x)\in Q,\\
        \partial_{\epsilon_i}u\lvert_{\epsilon=0}(0,x) = 0,\quad x\in\Omega,\\\partial_{\epsilon_i}u\lvert_{\epsilon=0}(t,x) = f_i(t,x), \quad (t,x)\in \Sigma.
    \end{cases}
\end{align*}
Now if we denote  $v_i(t,x):=\partial_{\epsilon_i}u\lvert_{\epsilon=0}$ ($1\leq i\leq m$), then the above IBVP reduces to
\begin{align}\label{IBVP when k=1}
    \begin{cases}
        \mathrm{i}\partial_t v_i(t,x) +\sum_{j=1}^{n}\left(\partial_j+\mathrm{i}A_j\right)^2 v_i(t,x)  + q(t,x)v_i(t,x)=0,~(t,x)\in Q,\\
        v_i(0,x) = 0,\quad x\in\Omega,\\v_i(t,x) = f_i(t,x), \quad (t,x)\in \Sigma.
    \end{cases}
\end{align}
From (Lemma 4.1, \cite{lai2024partial}), we have that $v_i \in H^{2\kappa}(Q)$ for $1\leq i\leq m$. The first order linearization of the DN map in \eqref{DN map for gover eqn 1}  is given by
\[\partial_{\epsilon_i}\Lambda_{A,q,B}(\epsilon f)\lvert_{\epsilon=0}= \left(\partial_{\epsilon_i} \partial_{\nu} u\lvert_{\Sigma}\right)\lvert_{\epsilon=0}=  \partial_{\nu}v_i\lvert_{\Sigma},\ 1\leq i\leq m.\]
Thus, it follows from \eqref{DN map for gover eqn 1} that 
\begin{align}\label{DN map when k=1}
    \partial_{\nu}v_i\lvert_{\Sigma}=\text{Known},\ \ \mbox{for any $f_i\in \widetilde{\mathcal{H}}^{2\kappa+\frac{3}{2}}_0(\Sigma)$}\ \mbox{and}\ 1\leq i\leq m.
\end{align}
Now since the Dirichlet data $f_i$ ($1\leq i\leq m$) in \eqref{IBVP when k=1} are part of measurement data,  therefore we denote 
$v_i:=v$  and $f_i:=\mathfrak{f}$, in \eqref{IBVP when k=1} and end up with the following    IBVP for $v$ 
\begin{align}\label{IBVP v}
    \begin{cases}
i\partial_t v(t,x) +\sum_{j=1}^{n}\left(\partial_j+\mathrm{i}A_j(t,x)\right)^2 v(t,x)  + q(t,x)v(t,x)=0,~(t,x)\in Q, \\
v(0,x) = 0,\quad x\in\Omega,\\
v(t,x) = \mathfrak{f}, \quad (t,x)\in \Sigma,
\end{cases}
\end{align}
and from first-order linearization of DN map $\Lambda_{A,q,B}$, we have
\begin{align}\label{DN map use for reconstructing a and q}
\Lambda_{A,q}(\mathfrak{f}):= \partial_{\nu} v \lvert_{\Sigma}=\text{Known}, \ \mbox{for any $\mathfrak{f}\in \widetilde{\mathcal{H}}^{2\kappa+\frac{3}{2}}_0(\Sigma)$}.
\end{align}
The equation \eqref{DN map use for reconstructing a and q} represent the DN map for the linearized IBVP given by \eqref{IBVP v}.

\subsubsection{\texorpdfstring{Reconstruction of \(A\)}{Reconstruction of A}}
To reconstruct $A$, we multiply \eqref{IBVP v} by $\overline{w}(t,x)$, where $w(t,x)$ satisfies the following backward problem  
     \begin{align}\label{adjoint equation 1}
         \mathrm{i}\partial_t w(t,x)+ \Delta w(t,x)=0, \ (t,x)\in Q  \ \mbox{and}\   w(T,x)=0, \ x\in\Omega
     \end{align}
    and integrate over $Q$, to yield 
     \begin{align*}
         \int_{Q} \left(\mathrm{i}\partial_{t} v +\Delta v +2\mathrm{i}A \cdot \nabla v-\lvert A\rvert^2 v + q v\right)(t,x) \overline{w}(t,x) \;dxdt=0.
     \end{align*}
  Next, the use of  integration-by-parts formula leads to 
  \begin{align*}
      \int_{\Omega} \mathrm{i}[v(T,x)\overline{w(T,x)}- v(0,x)\overline{w(0,x)}]\;dx-\int_{Q} \mathrm{i}v\partial_t \overline{w} \;dx dt +\int_{\Sigma}\overline{w}\partial_{\nu} v\; dS_xdt-\int_{\Sigma}v\partial_{\nu} \overline{w}\; dS_xdt \\+ \int_{Q} v\Delta \overline{w} \;dxdt+\int_{\Sigma} 2 \mathrm{i}(\nu \cdot A ) v \overline{w} ~dS_x dt    -2\mathrm{i}\int_{Q}(A\cdot \nabla \overline{w})v\;dx dt + \int_{Q} q_1 v \overline{w}\;dx dt=0
  \end{align*}
where $q_1:=q-\lvert A\rvert^2$. Using DN map given by \eqref{DN map use for reconstructing a and q} along with equations  \eqref{IBVP v} and \eqref{adjoint equation 1}, the above expression boils down to the following integral identity
\begin{align}\label{integral identity in case of linear}
    -2\mathrm{i}\int_{Q}(A\cdot \nabla \overline{w})v \;dx dt+ \int_{Q} q_1 v \overline{w}\;dx dt =\text{Known},
\end{align}
for all $v$ and $w$ solutions to \eqref{IBVP v} and   \eqref{adjoint equation 1}, respectively.
Next, to reconstruct the coefficient $A$, we use \textit{Geometric optics} (GO) solutions from (Section 3, \cite{kian2019holder}). 
 For $\varrho>0$, let $\Xi \in C_c^{\infty}(\mathbb{R},[0,1])$ with $\Xi=1$,  in $[2\varrho,T-2\varrho]$ and $\supp (\Xi)\subset (\varrho,T-\varrho)$ be smooth cut-off function satisfying 
\[\lVert \Xi \rVert_{W^{l,\infty}(\mathbb{R})} \lesssim  \varrho^{-l}\]
for each  $l\in \mathbb{N}$. Now for a fixed vector  $\o\in \mathbb{S}^{n-1}$ and $\lambda>0$ a large parameter, we define
\begin{align}\label{definition of u0 and u00}
\begin{aligned}
    &\phi(t,x):=\lambda(x \cdot \o-\lambda \lvert \o\rvert^2t),\ (t,x)\in Q,\\&T_v(t,x):=\Xi(t) \dfrac{\xi}{\lvert \xi\rvert}\cdot\nabla\left( e^{-\mathrm{i}(\tau t+ x \cdot \xi)} e^{-\mathrm{i}\int_{\mathbb{R}} \o\cdot A(t,x+s\o)~ds}\right)e^{\mathrm{i}\int_{0}^{\infty}\o\cdot A(t,x+s\o)ds}, \ (\tau,\xi) \in \mathbb{R}\times \o^{\perp} 
    \\&  T_w(t,x):=\Xi(t)
    \end{aligned}
\end{align} 
where $\displaystyle \omega^{\perp}:=\left\{x\in\mathbb{R}^n:\ x\cdot \omega=0\right\}$. 
Now following (Section 3, \cite{kian2019holder}), we choose  the GO solutions for  $v$ and $w$ of \eqref{IBVP v} and \eqref{adjoint equation 1} respectively, of the following form
\begin{align}\label{vi cgo 1}
v(t,x)=e^{\mathrm{i}\phi(t,x)} T_v(t,x)+R_v(t,x)
    \text{ and } w(t,x)=e^{\mathrm{i}\phi(t,x)}T_w(t,x)+R_w(t,x)
\end{align}
where
$v(0,\cdot)=w(T,\cdot)=0$, in $\O$. Also, the correction terms,  
$R_v$ and $R_w$ satisfy the following estimates
\begin{align}\label{reminder terms bound 1}
 \lambda\lVert R_v\rVert_{L^2(Q)}+\lVert \nabla R_v\rVert_{L^2(Q)}\leq C\quad \text{ and }\quad  \lambda\lVert R_w\rVert_{L^2(Q)}+\lVert \nabla R_w\rVert_{L^2(Q)}\leq C. 
\end{align}
Substitute \eqref{vi cgo 1} into the integral identity \eqref{integral identity in case of linear}, to arrive at  
\begin{align*}
   -2 \lambda \int_{Q}(\o\cdot A)T_v \overline{T_w}~dx dt  -2 \mathrm{i}\int_{Q}  e^{\mathrm{i}\phi}T_v(A\cdot \nabla\overline{R_w})  ~ dxdt-2 \mathrm{i}\int_{Q} R_v(A\cdot \nabla \overline{R_w})~dxdt
    \\-2 \lambda \int_{Q}(\o\cdot A)\overline{T_w}\left( e^{-\mathrm{i}\phi} R_v \right)~dxdt+\int_Q q_1 T_v\overline{T_w}(t,x)~dxdt\\+\int_Q q_1\left(e^{\mathrm{i}\phi} \overline{R_w}T_v+e^{-\mathrm{i}\phi} \overline{R_v}\overline{T_w}+R_v\overline{R_w}
\right)~ dxdt=\text{ Known}\end{align*} for any choice of $T_v$ and $T_w$ given by \eqref{definition of u0 and u00}. 
Multiply by $\lambda^{-1}$ in above expression and letting $\lambda \rightarrow\infty$ along with using equations \eqref{definition of u0 and u00} and \eqref{reminder terms bound 1}, we obtain
 \begin{align*}
\int_{Q} (\omega \cdot A)\Xi^2(t)\dfrac{\xi}{\lvert \xi\rvert}\cdot\nabla\left(  e^{-\mathrm{i}(\tau t+ x \cdot \xi)} e^{-\mathrm{i}\int_{\mathbb{R}} \o\cdot A(t,x+s\o)~ds}\right)e^{\mathrm{i}\int_{0}^{\infty}\o\cdot A(t,x+s\o)ds}~dxdt=  \text{ Known}
 \end{align*}
 for any $\Xi \in C_c^{\infty}(\mathbb{R},[0,1])$, $\xi\in\o^{\perp}$ and $\tau\in \mathbb{R}$. 
 Since $\Xi \in C_c^{\infty}(\mathbb{R},[0,1])$ and $A \in C_c^{\infty}(Q)$, therefore the above integral can be written as 
  \begin{align*}
\int_{\mathbb{R}^{1+n}} (\omega \cdot A)\Xi^2(t)\dfrac{\xi}{\lvert \xi\rvert}\cdot\nabla\left(  e^{-\mathrm{i}(\tau t+ x \cdot \xi)} e^{-\mathrm{i}\int_{\mathbb{R}} \o\cdot A(t,x+s\o)~ds}\right)e^{\mathrm{i}\int_{0}^{\infty}\o\cdot A(t,x+s\o)ds}~dxdt=  \text{ Known}
 \end{align*}
 for any $\Xi \in C_c^{\infty}(\mathbb{R})$, $\xi\in\o^{\perp}$ and $\tau\in \mathbb{R}$. 
 The use of decomposition $\mathbb{R}^n:= \mathbb{R}\omega \oplus \omega^{\perp}$ along with  the fact that $\Xi \in C_c^{\infty}(\mathbb{R},[0,1])$ and $\tau\in \mathbb{R}$ are arbitrary,  in the above integral, yields 
 \begin{align*}
  \int_{\omega^{\perp}} 
  \int_{\mathbb{R}} \left(\omega\cdot A(t,l+\sigma \omega)\right)\dfrac{\xi}{\lvert \xi\rvert} \cdot \nabla(e^{-\mathrm{i} \xi\cdot l} e^{-\mathrm{i}\int_{\mathbb{R}} \o\cdot A(t,l+s\o)~ds} )e^{ \mathrm{i}\int_{0}^{\infty} \omega\cdot A(t,l+\sigma \omega+s\omega) ds})d\sigma dl = \text{Known}
\end{align*}
for any $t\in [0,T]$ and $\xi\in \o^{\perp}$,  
where $dl$ denotes the Lebesgue measure on $\omega^{\perp}$. Now using the  change of variable $\sigma+s:= r$, in the above equation, we obtain that
\begin{align}\label{on the left side}
    \int_{\omega^{\perp}} 
    \int_{\mathbb{R}} e^{ \mathrm{i}\int_{\sigma}^{\infty} \omega\cdot A(t,l+r \omega) dr}\omega\cdot A(t,l+\sigma \omega)\dfrac{\xi}{\lvert \xi\rvert} \cdot \nabla e^{-\mathrm{i}( \xi\cdot l+\mathrm{i}\int_{\mathbb{R}} \o\cdot A(t,l+s\o)~ds}d\sigma dl= \text{Known}
\end{align}
for any $t\in [0,T]$ and $\xi\in \o^{\perp}$. Now if we denote 
\(\Psi:=\exp{\bigg( \mathrm{i}\int_{\sigma}^{\infty} \omega\cdot A(t,l+r\omega) dr}\bigg)\), 
then $\partial_{\sigma}\Psi= -\mathrm{i}\o \cdot A(t,l+\sigma \o) \Psi$. Integrating  this over $\mathbb{R}$ w.r.t. $\sigma$, yields that
\begin{align*}
-\mathrm{i}\int_{\mathbb{R}} \omega\cdot A(t,l+\sigma \omega) \Psi d\sigma=\int_{\mathbb{R}} \partial_{\sigma}\Psi \ d\sigma=1-  \exp{\bigg( \mathrm{i}\int_{\mathbb{R}} \omega\cdot A(t,l+r\omega) dr}\bigg).
\end{align*}
Substitute the above expression in \eqref{on the left side}, we get 
\begin{align}\label{2.12}
        \mathrm{i} \int_{\omega^{\perp}} \left(1-  e^{\mathrm{i}\int_{\mathbb{R}} \omega\cdot A(t,l+r\omega) dr}\right) \dfrac{\xi}{\lvert \xi\rvert} \cdot \nabla\left(e^{-\mathrm{i} \xi\cdot l} e^{-\mathrm{i}\int_{\mathbb{R}} \o\cdot A(t,l+s\o)~ds} \right)dl = \text{Known}
\end{align}
for any $t\in [0,T]$ and $\xi\in \o^{\perp}$. 
 Now since $A$ is compactly supported in $\Omega$ for each $t\in [0,T]$ and $\O\subset\mathbb{R}^n$ is bounded, therefore there exist $R>0$ such that  $\supp A(t,\cdot)\subset\O \subset B(0,R)$ and  
 \begin{align}\label{support condition}
     \supp\left(1-  e^{ \mathrm{i}\int_{\mathbb{R}}\omega\cdot A(t,l+r\omega) dr}\right) \subset \omega^{\perp}\cap B(0,R).
 \end{align}
The use of decomposition $\nabla g:=\nabla_{\perp} g+\o(\o\cdot \nabla g)$ and the integration by parts formula along with support condition \eqref{support condition} in \eqref{2.12} leads to 
\begin{align*}
 \text{Known}&= \mathrm{i} \int_{\omega^{\perp}} \left(1-  e^{\mathrm{i}\int_{\mathbb{R}} \omega\cdot A(t,l+r\omega) dr}\right) \dfrac{\xi}{\lvert \xi\rvert} \cdot \nabla\left(e^{-\mathrm{i} \xi\cdot l} e^{-\mathrm{i}\int_{\mathbb{R}} \o\cdot A(t,l+s\o)~ds} \right)dl \\&=\mathrm{i}
     \int_{\omega^{\perp}} \left(1-  e^{\mathrm{i}\int_{\mathbb{R}} \omega\cdot A(t,l+r\omega) dr}\right) \dfrac{\xi}{\lvert \xi\rvert} \cdot \nabla_{\perp}\left(e^{-\mathrm{i} \xi\cdot l} e^{-\mathrm{i}\int_{\mathbb{R}} \o\cdot A(t,l+s\o)~ds} \right)~dl\\&
     =\mathrm{i}
     \int_{\omega^{\perp}\cap B(0,R)} \left(1-  e^{\mathrm{i}\int_{\mathbb{R}} \omega\cdot A(t,l+r\omega) dr}\right) \dfrac{1}{\lvert \xi\rvert} \nabla_{\perp}\cdot \left(\xi e^{-\mathrm{i} \xi\cdot l} e^{-\mathrm{i}\int_{\mathbb{R}} \o\cdot A(t,l+s\o)~ds} \right)~dl
     \\&=-\mathrm{i}
     \int_{\omega^{\perp}} \left(e^{-\mathrm{i} \xi\cdot l} e^{-\mathrm{i}\int_{\mathbb{R}} \o\cdot A(t,l+s\o)~ds} \right)\dfrac{\xi}{\lvert \xi\rvert} \cdot \nabla_{\perp} \left(1-  e^{\mathrm{i}\int_{\mathbb{R}} \omega\cdot A(t,l+r\omega) dr}\right)~dl
   \\&
     =-\int_{\omega^{\perp}}e^{-\mathrm{i} \xi\cdot l}\dfrac{\xi}{\lvert \xi\rvert} \cdot \nabla_{\perp}\left(\int_{\mathbb{R}} \o\cdot A(t,l+s\o)~ds\right)~dl\\&=\int_{\omega^{\perp}}\int_{\mathbb{R}}e^{-\mathrm{i}\xi\cdot (l+s\o)}\dfrac{\xi}{\lvert \xi\rvert}\cdot \nabla_{\perp}\left( \o\cdot A(t,l+s\o)\right)~ds~dl\\&=\int_{\omega^{\perp}}\int_{\mathbb{R}}e^{-\mathrm{i}\xi\cdot (l+s\o)}\dfrac{\xi}{\lvert \xi\rvert}\cdot \nabla\left( \o\cdot A(t,l+s\o)\right)~ds~dl=
   \int_{\mathbb{R}^n}e^{-\mathrm{i}\xi\cdot x}\dfrac{\xi}{\lvert \xi\rvert}\cdot \nabla\left( \o\cdot A(t,x)\right)~dx
\end{align*}
for any $t\in [0,T]$ and $\xi\in \o^{\perp}$.
Using Fourier transform of the space variable, the above expression reduces to
\begin{align}\label{omega dot A known}
\text{Known}=\int_{\mathbb{R}^{n}}e^{-\mathrm{i}\xi\cdot x}\dfrac{\xi}{\lvert \xi\rvert} \cdot \nabla \left( \o\cdot A(t,x)\right)dx =\mathrm{i}\lvert \xi\rvert \o\cdot \widehat{A}(t,\xi)
\end{align}
for any $t\in [0,T]$ and $\xi\in \o^{\perp}$.
Also, using the Fourier transform along with the given hypothesis $\nabla \cdot A=0$, we have
$ \xi \cdot   \widehat{A}(t,\xi)=0$, for all $\xi \in \mathbb{R}^n\setminus\{0\}$.  Also for any $\xi\in \mathbb{R}^n$, we choose $\o_i\in  \mathbb{S}^{n-1}$ for $1\leq i\leq n-1$ such that 
 $\left\{\o_1,\cdots,\o_{n-1},\frac{\xi}{\lvert \xi\rvert}\right\}$ form an orthonormal basis for $\mathbb{R}^n$. Using these basis elements,
 we can write $\widehat{A}(t,\xi)=\sum_{i=1}^{n-1}a_i \o_i+a_n\frac{\xi}{\lvert \xi\rvert}$ where $a_i:=\langle \widehat{A}(t,\xi),\o_i\rangle$ for $1\leq i\leq n$ and $a_n=\langle \widehat{A}(t,\xi),\frac{\xi}{\lvert \xi\rvert}\rangle$.
 Since $ \xi \cdot   \widehat{A}(t,\xi)=0$, invoking \eqref{omega dot A known}, we have 
 \[
\widehat{A}(t,\xi)=\sum_{i=1}^{n-1}
 \langle \widehat{A}(t,\xi),\o_i\rangle \o_i=\text{ Known}, \ \text{for any $t\in [0,T]$ and $\xi \in \mathbb{R}^n$.}
 \]
  Now, by after utilizing the Fourier inversion formula, we get 
\(A(t,x)=\text{Known}, \ \text{for all }(t,x) \in Q.\)\qed

\subsubsection{\texorpdfstring{Reconstruction of \(q\)}{Reconstruction of q}} To reconstruct $q$, we  multiply equation \eqref{IBVP v} by $\overline{w}(t,x)$, where $w(t,x)$ satisfies the following backward problem
     \begin{align}\label{modified adjoint equation 2}
     \begin{aligned}
         \mathrm{i}\partial_t w + \Delta w+2\mathrm{i} A\cdot \nabla w=0, \ (t,x)\in Q \ \mbox{and}\  w(T,x)=0, \ x\in\Omega
         \end{aligned}
     \end{align}
    and integrating over $Q$, we get
   \begin{align*}
         \int_{Q} \left(\mathrm{i}\partial_{t} v(t,x) +\Delta v(t,x) +2\mathrm{i}A \cdot \nabla v(t,x) +q_1(t,x) v(t,x)\right) \overline{w(t,x)} \;dxdt=0
     \end{align*}
 where $q_1=q-\lvert A \rvert^2$.  Next, the use of integration-by-parts formula in the above equation, leads to 
  \begin{align*}
      \int_{\Omega} \mathrm{i}\left[v(T,x)\overline{w}(T,x)- v(0,x)\overline{w}(0,x)\right]\;dx-\int_{Q} \mathrm{i}v\partial_t \overline{w} \;dx dt +\int_{\Sigma}\overline{w}\partial_{\nu} v\; dS_xdt-\int_{\Sigma}v\partial_{\nu} \overline{w}\; dS_xdt \\+ \int_{Q} v\Delta \overline{w} \;dxdt+\int_{\Sigma} 2 \mathrm{i}(\nu \cdot A ) v \overline{w} ~dS_x dt    -2\mathrm{i}\int_{Q}(A\cdot \nabla \overline{w})v\;dx dt + \int_{Q} q_1 v \overline{w}\;dx dt=0.
  \end{align*}
Using DN map \eqref{DN map use for reconstructing a and q} along with equations \eqref{IBVP v} and \eqref{modified adjoint equation 2}, the above expression simplifies to 
\begin{align}\label{q_1 integral}
       \int_{Q}q_1(t,x) v(t,x) \overline{w(t,x)}\;dxdt= \text{Known}
\end{align}
for any $v$ and $w$ solutions to \eqref{IBVP v} and \eqref{modified adjoint equation 2}, respectively. 
Now as before, we choose the GO solutions for $v$ and $w$, solutions to equations  \eqref{IBVP v} and \eqref{modified adjoint equation 2} respectively and taking the  following form
\begin{align}
v(t,x)=e^{\mathrm{i}\phi(t,x)} \Xi(t)e^{-\mathrm{i}(\tau t+ x\cdot \xi)}e^{\mathrm{i}\int_{0}^{\infty}\o\cdot A(t,x+s\o)ds}+R_v(t,x)
\end{align}
and
\begin{align}
w(t,x)=e^{\mathrm{i}\phi(t,x)}\Xi(t)e^{\mathrm{i}\int_{0}^{\infty}\o\cdot A(t,x+s\o)ds}+R_w(t,x)
\end{align}
where
$v(0,\cdot)=w(T,\cdot)=0$ in $\O$. Also, the correction terms,  
$R_v$ and $R_w$ satisfy the following estimates
\begin{align}\label{estimate on remainder terms}
 \lambda\lVert R_v\rVert_{L^2(Q)}+\lVert \nabla R_v\rVert_{L^2(Q)}\leq C\quad \text{ and }\quad  \lambda\lVert R_w\rVert_{L^2(Q)}+\lVert \nabla R_w\rVert_{L^2(Q)}\leq C. 
\end{align}
Using aforementioned $v$ and $w$ in  \eqref{q_1 integral}, we get
\begin{align}
    \int_{Q} q_1 \Xi^2(t)&  e^{-\mathrm{i}(\tau t + x \cdot \xi)}\ dx dt+\int_{Q} q_1 \Xi(t)R_v e^{-\mathrm{i}\phi}e^{-\mathrm{i}\int_{0}^{\infty}\o\cdot A(t,x+s\o)ds}\ dx dt\\&+\int_{Q} q_1 \Xi(t) \overline{R_w}e^{\mathrm{i}\phi}e^{-\mathrm{i}(\tau t + x \cdot \xi)}e^{\mathrm{i}\int_{0}^{\infty}\o\cdot A(t,x+s\o)ds}\ dx dt+   \int_{Q} q_1 R_v\overline{R_w}\ dx dt=\text{ Known}.
\end{align}
After taking $\lambda \rightarrow\infty$ along with using equation \eqref{estimate on remainder terms}, we get 
\begin{align}
   \int_{Q}\Xi^2(t) q_1(t,x) e^{-\mathrm{i}(\tau t + x \cdot \xi)}\ dx dt= \text{ Known}, 
   \ \text{for any $\Xi \in C_c^{\infty}(\mathbb{R},[0,1])$, $\xi\in\omega^{\perp}$ and $\tau\in \mathbb{R}$}.
\end{align}
Next  utilize the fact that  $\Xi \in C_c^{\infty}(\mathbb{R},[0,1])$ is arbitrary  and $q_1 \in C_c^{\infty}(Q)$, in the  above integral, to obtain
\[
\int_{\mathbb{R}^{n}}  q_1(t,x) e^{-\mathrm{i} x \cdot \xi}\ dx=\text{Known},\ \text{for any $t\in [0,T]$ and $\xi\in \omega^{\perp}$}.
\]
Now using the fact that $\omega\in \mathbb{S}^{n-1}$ is arbitrary and the fact that $\displaystyle \int_{\mathbb{R}^{n}}  q_1(t,x) e^{-\mathrm{i} x \cdot \xi}\ dx$ is known for any $\xi\in \omega^{\perp}$ and $t\in [0,T]$,  we obtain that\[\int_{\mathbb{R}^{n}}  q_1(t,x) e^{-\mathrm{i} x \cdot \xi}\ dx=\text{known}, \ \text{for any $\xi\in \mathbb{R}^n$ and $t\in [0,T]$}.\] 
Finally, using the Fourier inversion formula, we get 
\[q_1(t,x)=\text{Known}, \  \text{ for all }(t,x) \in Q.\] 
Since $A$ is already been reconstructed therefore we conclude that $q(t,x)$ is known for all $(t,x)\in Q$.\qed
\subsection{Reconstruction of nonlinearity coefficients}\label{sub2} This subsection is devoted to deriving the reconstruction formulae for $r_i$ $(1\leq i\leq m)$. As mentioned earlier, following 
\cite{lassas2021inverse, kang2002identification, lassas2020partial, nakamura2008inverse, bhardwaj2026reconstruction}, we use the higher-order linearization of solution to the nonlinear equation  for establishing the reconstruction of  these coefficients.

\subsubsection{Higher order linearization} The $k$-th ($2\leq k\leq m$) order linearization is done by applying $\partial_{\epsilon_{j_1}\cdots\epsilon_{j_k}}^k$  ($j_1,\cdots,j_k\in\{1,2,\cdots,m\}$) to  \eqref{govern pde with boundary data ef}  and substituting $\epsilon=0$. 
Hence, we start by applying the partial differential operator $\partial_{\epsilon_1\cdots \epsilon_k}^k$ $(2\leq k\leq m)$  to  \eqref{govern pde with boundary data ef} and  obtain that  
\begin{align}\label{upto m deivative of gover eqn}
\begin{aligned}
    \begin{cases}
     \mathrm{i}\partial_t\left( \partial^k_{\epsilon_1\cdots \epsilon_k}u\right) +\sum_{j=1}^{n}\left(\partial_j+\mathrm{i}A_j\right)^2 \left( \partial^k_{\epsilon_1\cdots \epsilon_k}u\right)  + q\left( \partial^k_{\epsilon_1\cdots \epsilon_k}u\right)=\partial^k_{\epsilon_1\cdots \epsilon_k} B(t,x,u,\overline{u}),~(t,x)\in Q,\\   
\partial^k_{\epsilon_1\cdots \epsilon_k}u(0,x) = 0,\quad x\in\Omega,\\ \partial^k_{\epsilon_1\cdots \epsilon_k}u(t,x) = 0, \quad (t,x)\in \Sigma.
    \end{cases}
    \end{aligned} 
\end{align} 
Now using well-posedness of the IBVP \eqref{govern pde with boundary data ef}, we get  that 
\(\partial^k_{\epsilon_1\cdots \epsilon_k}\rvert_{\epsilon=0} \left(B(t,x,u,\overline{u})\right)=0, \) for any $2\leq k\leq m-1$. Thus, after evaluating
\eqref{upto m deivative of gover eqn} at $\epsilon=0$, along with denoting by $z_k:=\partial^k_{\epsilon_1\cdots \epsilon_k}u\lvert_{\epsilon=0}$, for $2\leq k\leq m-1$, we obtain the following IBVP for $z_k$ ($2\leq k\leq m-1$)
\begin{align}\label{intermediate pde}
    \begin{cases}
     \mathrm{i}\partial_t z_k(t,x) +\sum_{j=1}^{n}\left(\partial_j+\mathrm{i}A_j(t,x)\right)^2 z_k(t,x) + q(t,x)z_k(t,x)=0,~(t,x)\in Q,\\   
     z_k(0,x) = 0,\quad x\in\Omega,\\z_k(t,x) = 0, \quad (t,x)\in \Sigma.
    \end{cases}
\end{align}
 The well-posedness of the IBVP \eqref{intermediate pde} implies that the $z_k\equiv 0$, in $Q$ for any $2\leq k\leq m-1$. Finally, applying $\displaystyle \partial^{m}_{ \epsilon_1\epsilon_2\cdots\epsilon_m}$ to \eqref{govern pde with boundary data ef} and  evaluating at $\epsilon = 0$, we obtain
  \begin{align}\label{IBVP when k=m}
\begin{cases}
\mathrm{i}\partial_t X +\sum_{j=1}^{n}\left(\partial_j+\mathrm{i}A_j\right)^2 X + qX= \partial^m_{\epsilon_1\cdots \epsilon_m}\rvert_{\epsilon=0} \left(B(t,x,u,\overline{u})\right):=B_1(t,x),~(t,x)\in Q, \\
X(0,x) = 0,\quad x\in\Omega,\\
X(t,x) = 0, \quad (t,x)\in \Sigma,
\end{cases}
\end{align}
where $X:=\partial^m_{\epsilon_1\cdots \epsilon_m}u\lvert_{\epsilon=0}$ and the function  $B_1(t,x)$ is given by
\begin{align*}
 B_1(t,x):&=\partial^m_{\epsilon_1\cdots \epsilon_m}\rvert_{\epsilon=0} \left(B(t,x,u,\overline{u})\right)=\partial^m_{\epsilon_1\cdots \epsilon_m}\rvert_{\epsilon=0}\left[\sum_{a=1}^{m}r_a(t,x)u^{a}(t,x)(\overline{u})^{m-a}(t,x)\right]\\&=m! r_m\prod_{k=1}^m v_k+\sum_{a=1}^{m-1}a!(m-a)!r_a\sum_{\substack{I\subset\{1,2,\cdots,m\}
\\
\lvert I \rvert =a}}\left(\prod_{j \in I}v_j\right)\left(\prod_{j \notin I}\overline{v_j}\right)
\end{align*}
where $v_k$ for $1\leq k\leq m$ are solutions to IBVP \eqref{IBVP when k=1}. The $m$-th order linearization of the DN map \eqref{DN map for gover eqn 1} is given by
\[\partial^m_{\epsilon_1\cdots \epsilon_m}\Lambda_{A,q,B}(\epsilon f)\lvert_{\epsilon=0}= \left(\partial^m_{\epsilon_1\cdots \epsilon_m} \partial_{\nu} u\lvert_{\Sigma}\right)\lvert_{\epsilon=0}=  \partial_{\nu}X\lvert_{\Sigma}  \ \mbox{ for any } \ \epsilon f\in \mathcal{E}_{\delta}(\Sigma).\]
Thus, it follows from \eqref{DN map for gover eqn 1} that 
\begin{align}\label{DN map when k=m}
    \partial_{\nu}X\lvert_{\Sigma}=\text{Known}.
\end{align}
With this information, we are now in a position to establish the reconstruction formulae for the nonlinearity coefficients $r_a$ for $1\leq a\leq m$. We split the proof for reconstruction of these coefficients into two parts: the first focuses on reconstructing $r_m$, while the second addresses the reconstruction of $r_a$ for $1\leq a\leq m-1$.

\subsubsection{\texorpdfstring{Reconstruction of \(r_m\)}{Reconstruction of rm}} To reconstruct $r_m$, we multiply equation \eqref{IBVP when k=m} by $\overline{w}(t,x)$, where $w(t,x)$ satisfies the following backward problem
     \begin{equation}\label{adjoint equation v power m plus one}
       \mathcal{P}_{A,q}  w(t,x)=0,~(t,x)\in Q, \ \mbox{and}\ w(T,x)=0, \ x \in \O
     \end{equation}
    and integrating over $Q$, we have
    \[\int_{Q} \bigg(\mathrm{i}\partial_t X +\Delta X+2\mathrm{i}A\cdot \nabla X-\lvert A\rvert^2 X  + qX -B_1 \bigg)(t,x) \overline{w}(t,x)~ dx dt =0.\]
After using the integration-by-parts formula in above equation, we obtain
  \begin{align*}
      &\int_{\Omega} \mathrm{i}\left[X(T,x)\overline{w}(T,x)- X(0,x)\overline{w}(0,x)\right]\;dx-\int_{Q} \mathrm{i}X(t,x)\partial_t \overline{w}(t,x) \;dx dt\\
       & \quad +\int_{\Sigma}\left[\overline{w}(t,x)\partial_{\nu} X(t,x)-X(t,x)\partial_{\nu} \overline{w}(t,x)\right] dS_xdt + \int_{Q} X(t,x)\Delta \overline{w}(t,x) \;dxdt \\
     & \quad +\int_{\Sigma} 2 \mathrm{i}(\nu \cdot A )(t,x) X(t,x) \overline{w}(t,x) ~dS_x dt    -2\mathrm{i}\int_{Q}(A\cdot \nabla \overline{w})(t,x)X(t,x)\;dx dt \\ 
     & \quad - \int_{Q}\left[ \lvert A\rvert^2 +q\right](t,x) X(t,x)\overline{w}(t,x)\;dx dt -\int_{Q} B_1(t,x) \overline{w}(t,x)~ dx dt=0.
  \end{align*}
Now utilizing equations \eqref{DN map when k=m}, \eqref{adjoint equation v power m plus one}, and \eqref{IBVP when k=m}, the foregoing expression simplifies to the following integral identity.
\begin{align}\label{integral identity}
\begin{aligned}
      &\int_{Q} B_1(t,x)\overline{w}(t,x)\;dxdt=m!\int_{Q}r_m(t,x)\prod_{k=1}^m v_k(t,x)dxdt\\ & \quad + \sum_{a=1}^{m-1}a!(m-a)!\int_{Q}r_a(t,x)\sum_{\substack{I\subset\{1,2,\cdots,m\}\\
\lvert  I \rvert =a}}\left(\prod_{j \in I}v_j\right)\left(\prod_{j \notin I}\overline{v_j}\right) \overline{w}(t,x) \;dxdt=\text{known} 
\end{aligned}
\end{align}
for any $v_k\ (1\leq k\leq m)$ and $w$ solutions to \eqref{IBVP when k=1} and \eqref{adjoint equation v power m plus one}, respectively. 
From the above integral identity, we reconstruct the coefficients $r_a$ for $1\leq a\leq m$. To this end, we use the GO solutions with \textit{concentrated amplitudes} and we refer to \cite{lai2024partial} for construction of such solutions.  
  We also refer the reader to  \cite{kian2019holder,  lassas2025coefficient, bhardwaj2026reconstruction}, in which similar approaches have been utilized.
We start with fixing  $\o_l \in \mathbb{R}^n\setminus\{\textbf{0}\}$ and  denote $\phi_l:=\lambda(x \cdot \o_l-\lambda \lvert \o_l\rvert^2t),$ for $0\leq l\leq m$, where $\lambda>0$ is a large parameter. 
Next for each $\o_{l}$ ($1\leq l\leq m$), we choose $\xi_1^l, \cdots\xi_{n-1}^l\in \mathbb{R}^n$ such that $\left\{\dfrac{\o_l}{\lvert{\o_l}\rvert},\xi_1^l, \cdots\xi_{n-1}^l\right\}$ forms an orthonormal basis for $\mathbb{R}^n$.
Also choose  $\mathcal{T} \in C^{\infty}_c(0,T)$ and $\chi \in C^{\infty}_c(\mathbb{R},[0,1])$,  the smooth cutoff functions such that $\chi(x)=1$, for $\lvert x\rvert\leq\frac{1}{2}$ and $\supp(\chi) \subset B(0,1)$.
Now for a fixed $x_0\in \O$ and $h>0$ sufficiently small, we define
\begin{align}\label{concentrated amplitude}
\begin{aligned}
   U_0^l(t,x):=\mathcal{T}(t)\prod_{j=1}^{n-1} \chi\left(
    \dfrac{(x-x_0)\cdot \xi_j^l}{h}
    \right) \exp\left({\mathrm{i} \int_{0}^{\infty} \o_l\cdot A(t,x+r \o_l) ~dr}\right),\ (t,x)\in Q
    \end{aligned}
\end{align}
for $0\leq l\leq m$.
Also, for fixed $\o\in \mathbb{S}^{n-1}$ and $x:=y+(x\cdot\o)\o \in \O$, we define
\begin{align}\label{concentrated amplitude 1}
    U_j^l(t,x)&:=\mathrm{i} \int_{-x\cdot\o}^0 \exp\left(-\mathrm{i}\int_{r}^0\o \cdot A(t,x+r_1\o)~dr_1\right)\mathcal{P}_{A,q} U_{j-1}^l (t,x+r\o)~ dr \ (1\leq j\leq N)
\end{align}
with initial condition $U_j^l(t,x)\big\lvert_{x\cdot\o=0}=0$ for all $t\in (0,T)$.   In the expression above, the linear operator $\mathcal{P}_{A,q}$ is as defined in \eqref{linear operator P}.
Now following (Proposition $3.2$, \cite{lai2024partial}), we choose  the GO solutions for  $v_k$ and $w$ solving \eqref{IBVP when k=1} and \eqref{adjoint equation v power m plus one} respectively, and taking the following form
\begin{align}\label{vi cgo}
v_k(t,x)=e^{\mathrm{i}\phi_k(t,x)}\left( U_0^k(t,x)+\sum_{j=1}^{N} \lambda^{-j}U_j^k (t,x)\right)+R_{v_k}(t,x),\ \text{for $1\leq k\leq m$}
\end{align}
and
\begin{align}\label{w cgo}
w(t,x)=e^{\mathrm{i}\phi_0(t,x)}\left( U_0^0(t,x)+\sum_{j=1}^{N} \lambda^{-j}U_j^0 (t,x)\right)+R_{w}(t,x)
\end{align}
where $v_k(0,\cdot)=0,$ for $1\leq k \leq m$ and $w(T,\cdot)=0$, in $\O$. For any natural number $s$ and $0\leq l\leq m$, we have
\begin{align}\label{reminder terms bound}
\begin{aligned}
    \lVert U_j^l \rVert_{H^s(Q)}&=\mathcal{O}(1), \quad 0\leq j\leq N , \\
    \lVert R_{v_k}\rVert_{H^s(Q)}=\mathcal{O}(\lambda^{-N+2s}) \ &\text{ and }  \lVert R_{w}\rVert_{H^s(Q)}=\mathcal{O}(\lambda^{-N+2s}).
    \end{aligned}
\end{align}
 Next, we substitute \eqref{vi cgo} and \eqref{w cgo} into the integral equation \eqref{integral identity} and use the Sobolev embedding theorem to arrive at 
\begin{align}\label{reduced integral identity}
\begin{aligned}
 \text{Known}&= \displaystyle \int_{Q}    m! r_me^{\mathrm{i}\left(\sum_{k=1}^m\phi_k-\phi_0\right)} \left(\prod_{k=1}^m U_0^k\right)\overline{U_0^0}~ dx dt\\&\quad+
 \displaystyle \sum_{a=1}^{m-1}\sum_{\substack{I\subset\{1,2,\cdots,m\}\\
\lvert  I \rvert =a}}\int_{Q}
       a!(m-a)!r_ae^{\mathrm{i}\left(\sum_{j\in I}\phi_j- \sum_{j\notin I}\phi_j-\phi_0\right)}\left(\prod_{j \in I}U_0^j\right)\left(\prod_{j \notin I}\overline{U_0^j}\right) \overline{U_0^0} \;dxdt\\&\quad+\mathcal{O}(\lambda^{-1})
\end{aligned}
\end{align}
for $s\geq  \kappa+1$ and $N> 2s$, by virtue of \eqref{reminder terms bound}. In the next step, we select $\phi_l$ for $0\leq l\leq m$ so that the exponential factor in the first term of \eqref{reduced integral identity} vanishes while  in the  second term of \eqref{reduced integral identity} it must remain non-vanishing for every admissible set $I$.
To reconstruct $r_m$, we choose $\o_k$ for $0\leq k\leq m$ so that 
\begin{align}\label{omega sum in reconstruction of rm}
    \o_1+\cdots+\o_{m}=\o_0\ \text{ and }\lvert \o_1\rvert^2+\cdots \lvert \o_{m}\rvert^2=\lvert \o_0\rvert^2.
\end{align}
\begin{table}[ht]
\centering
\begin{tabular}{|Sc|Sc|}
   \hline $m= 2$ & $m\geq3 $  \\\hline
     $  \o_1=(1,0,\cdots,0)$& $\o_1=\o_2=\cdots=\o_{m-1}=(0,-1,0,\cdots,0) $\\\hline$\o_2=(0,-1,0,\cdots,0)$&$\o_{m}=\left(0,\dfrac{m-2}{2},0,\cdots,0\right)$ \\\hline$\o_0=(1,-1,0\cdots,0)$ &$\o_0=\left(0,-\dfrac{m}{2},0\cdots,0\right)$\\\hline
\end{tabular}
\vspace{0.5cm}
\caption{A particular choice of $\o_k$ for $0\leq l\leq m$ to reconstruct $r_m$}
\label{table1}
\end{table}
A particular choice of $\omega_k$ ($0\leq k\leq m$) satisfying \eqref{omega sum in reconstruction of rm} is given in  \Cref{table1}, below.  Hence, building on this choice of $\omega_k's$,  we observe that
\begin{align}
    \phi_1+\cdots+\phi_m-\phi_0= \lambda x\cdot(\o_1+\cdots+\o_{m}-\o_0)- \lambda^2 t(\lvert \o_1\rvert^2+\cdots \lvert \o_{m}\rvert^2-\lvert \o_0\rvert^2)=0.
\end{align}
Furthermore, for these choices of $\o_i$, $\sum_{j\in I}\phi_j- \sum_{j\notin I}\phi_j-\phi_0\neq0$ for any $I\subset \{1,2,\cdots,m\}$ with  $\lvert I \rvert =a \ \text{for}\ 1\leq a\leq m-1$. If possible, let us suppose that $\sum_{j\in I}\phi_j- \sum_{j\notin I}\phi_j-\phi_0=0$, where $I=\{i_1,\cdots ,i_a\}$, then 
\begin{align}\label{permutation of phi i}
\Phi:=\phi_{i_1}+\cdots +\phi_{i_a}-(\phi_{j_1}+\cdots+\phi_{j_{m-a}}+\phi_0)=0
\end{align}
where $\{i_1,\cdots ,i_a,j_1,\cdots ,j_{m-a}\}$ is some permutation of $\{1,2,\cdots,m\}$. From definition of $\phi_l$ for $0\leq l\leq m$, the foregoing expression reduces to  
\begin{align}\label{omega sum in reconstruction of rm for |l|=a}
\begin{aligned}
   \omega_{i_1}+\cdots+  \omega_{i_a}&=\o_{j_1}+\cdots+\o_{j_{m-a}}+\o_0,\\
   \lvert  \omega_{i_1}\rvert^2+\cdots+\lvert  \omega_{i_a}\rvert^2&=\lvert  \omega_{j_1}\rvert^2+\cdots+\lvert  \omega_{j_{m-a}}\rvert^2+\lvert  \omega_{0}\rvert^2.
   \end{aligned}
\end{align}
Since $\{i_1,\cdots ,i_a,j_1,\cdots ,j_{m-a}\}$ is some permutation of $\{1,2,\cdots,m\}$, we can rewrite \eqref{omega sum in reconstruction of rm} as
\begin{align}\label{rewritten omega sum in reconstruction of rm}
    \begin{aligned}
   \omega_{i_1}+\cdots+  \omega_{i_a}&+\o_{j_1}+\cdots+\o_{j_{m-a}}=\o_0,\\
   \lvert  \omega_{i_1}\rvert^2+\cdots+\lvert  \omega_{i_a}\rvert^2&+\lvert  \omega_{j_1}\rvert^2+\cdots+\lvert  \omega_{j_{m-a}}\rvert^2=\lvert  \omega_{0}\rvert^2.
   \end{aligned}
\end{align}
The difference of \eqref{omega sum in reconstruction of rm for |l|=a} and \eqref{rewritten omega sum in reconstruction of rm} leads to
\[\o_{j_1}+\cdots+\o_{j_{m-a}}=0 \ \text{ and }\ \lvert \o_{j_1} \rvert^2+\cdots +\lvert \o_{j_{m-a}} \rvert^2 =0.\]
The aforementioned expression implies that $\o_{j_1}=\cdots=\o_{j_{m-a}}=0$, which is
 a contradiction. Thus, the exponential factor in the second term of \eqref{reduced integral identity} is non-vanishing for every admissible set $I$. 
 \begin{step}\label{step 4} To proceed further, we first verify the hypothesis needed to apply \textit{nonstationary phase lemma} for oscillating integral. Note that the phase function $\Phi$ satisfies
\begin{align*}
    \partial_t\Phi&=-\lambda^2\left(\lvert \o_{i_1}\rvert^2+\cdots+\lvert \o_{i_a}\rvert^2-(\lvert \o_{j_1}\rvert^2+\cdots+\lvert \o_{j_{m-a}}\rvert^2+\lvert \o_{0}\rvert^2)\right),\\ 
   \nabla \Phi&= \lambda(\o_{i_1}+\cdots+\o_{i_a}-(\o_{j_1}+\cdots \o_{j_{m-a}}+\o_0)).
\end{align*}
 As a result, $\nabla_{t,x}\Phi:=(\partial_t,\nabla)\Phi\neq(0,\textbf{0})$ for any $(t,x) \in Q$ and $\left\lvert\nabla_{t,x}\Phi\right\rvert=\mathcal{O}(\lambda^2).$
 Also, from the construction of $U_0^j~(0\leq j \leq m)$ and assumption on $r_a$, we have
 \begin{align*}
     r_a\left(\prod_{j \in I}U_0^j\right)\left(\prod_{j \notin I}\overline{U_0^j}\right) \overline{U_0^0} \in C^{\infty}_c(Q).
 \end{align*}
Following (Lemma 3.14, \cite{zworski2012semiclassical}), for any natural number $N$, there exist $C_{Q,N}>0,$ such that
\begin{align*}
  \left\lvert  \int_{Q}r_a e^{\mathrm{i}\Phi}\left(\prod_{j \in I}U_0^j\right)\left(\prod_{j \notin I}\overline{U_0^j}\right) \overline{U_0^0} \;dxdt\right\rvert\leq C_{Q,N}\lambda^{-2N}.
\end{align*}
Thus, we conclude that 
 \begin{align}\label{order of lambda minus}
    \left\lvert  \sum_{\substack{I\subset\{1,2,\cdots,m\}\\
\lvert  I \rvert =a}}\int_{Q}
       a!(m-a)!r_ae^{\mathrm{i}\left(\sum_{j\in I}\phi_j- \sum_{j\notin I}\phi_j-\phi_0\right)}\left(\prod_{j \in I}U_0^j\right)\left(\prod_{j \notin I}\overline{U_0^j}\right) \overline{U_0^0} \;dxdt\right\rvert
       =\mathcal{O}(\lambda^{-2}).
 \end{align}
Taking $\lambda \rightarrow\infty$ along with using above expression in \eqref{reduced integral identity}, to obtain that 
 \begin{align}
   \text{Known}&= \int_{Q}r_m(t,x)\left(\prod_{k=1}^m U_0^k\right)(t,x) \overline{U_0^0}(t,x) \;dxdt\\
  &=\int_{Q}r_m(t,x)\prod_{k=1}^m \left[ \mathcal{T}(t)\prod_{j=1}^{n-1} \chi\left(
    \dfrac{(x-x_0)\cdot \xi_j^k}{h}
    \right) \exp\left({\mathrm{i} \int_{0}^{\infty} \o_k\cdot A(t,x+r \o_k) ~dr}\right)
  \right]\times
  \\&\qquad \left[ \mathcal{T}(t)\prod_{j=1}^{n-1} \chi\left(
    \dfrac{(x-x_0)\cdot \xi_j^0}{h}
    \right) \exp\left({-\mathrm{i} \int_{0}^{\infty} \o_0\cdot A(t,x+r \o_0) ~dr}\right)\right]~dxdt 
\end{align}
Since $\mathcal{T}\in C_c^{\infty}(0,T)$ is arbitrary and $r_m \in C_c^{\infty}(Q)$, we have
 \begin{align}\label{r_m integral identity}
\begin{aligned}
\text{Known}&=\int_{\O}r_m(t,x)\left(\prod_{k=1}^m \left[\prod_{j=1}^{n-1} \chi\left(
    \dfrac{(x-x_0)\cdot \xi_j^k}{h}
    \right)\right] \prod_{j=1}^{n-1} \chi\left(
    \dfrac{(x-x_0)\cdot \xi_j^0}{h}
    \right)\right)
  \times
  \\&\qquad \left[  \exp\left({\mathrm{i}\sum_{k=1}^m  \int_{0}^{\infty} \o_k\cdot A(t,x+r \o_k) ~dr
  -\mathrm{i} \int_{0}^{\infty} \o_0\cdot A(t,x+r \o_0) ~dr}\right)\right]~dx\\&=\int_{B(x_0,h)}r_m(t,x)\left(\prod_{k=1}^m \left[\prod_{j=1}^{n-1} \chi\left(
    \dfrac{(x-x_0)\cdot \xi_j^k}{h}
    \right)\right] \prod_{j=1}^{n-1} \chi\left(
    \dfrac{(x-x_0)\cdot \xi_j^0}{h}
    \right)\right)
  \times
  \\&\qquad \left[  \exp\left({\mathrm{i}\sum_{k=1}^m  \int_{0}^{\infty} \o_k\cdot A(t,x+r \o_k) ~dr
  -\mathrm{i} \int_{0}^{\infty} \o_0\cdot A(t,x+r \o_0) ~dr}\right)\right]~dx
    \end{aligned}
\end{align}
for all $t\in (0,T)$, where in the last step of the above expression, we have used the support condition on $\chi$. 
 Now multiply the above equation by $h^{-n}$  and letting $h \rightarrow 0$, we obtain
\begin{align}
    \text{Known}&=r_m(t,x_0) \left[  \exp\left({\mathrm{i}\sum_{k=1}^m  \int_{0}^{\infty} \o_k\cdot A(t,x_0+r \o_k) ~dr
  -\mathrm{i} \int_{0}^{\infty} \o_0\cdot A(t,x_0+r \o_0) ~dr}\right)\right]
\end{align}
where to arrive at the above equation, we used the properties of  $\chi$. 
Finally, utilizing the fact that $A$ is known and the point $x_0\in\Omega$ is arbitrary, we get 
\begin{align}
    r_m(t,x)= \text{Known, for all $(t,x)\in Q$}.
\end{align} 
 \end{step}

\subsubsection{\texorpdfstring{Reconstruction of \(r_l~(1\leq l\leq m-1)\)}{Reconstruction of rl}} To reconstruct $r_l$ for $1\leq l\leq m-1$, we start with rewriting 
\eqref{integral identity} as below 
\begin{align}\label{r term}
\begin{aligned}
      \text{Known}&= \int_{Q}\left(
       m! r_m \prod_{k=1}^m v_k\right)\overline{w}~dxdt+\sum_{a=1}^{m-1}\sum_{\substack{I\subset\{1,2,\cdots,m\}\\
\lvert  I \rvert =a}}\displaystyle \int_{Q}
     a!(m-a)!r_a\left(\prod_{j \in I}v_j\right)\left(\prod_{j \notin I}\overline{v_j}\right) \overline{w} \;dxdt\\& =\int_{Q}\left(
       m! r_m \prod_{k=1}^m v_k\right)\overline{w}~dxdt+\int_{Q}
     l!(m-l)!r_l \left(\prod_{j=1}^l v_j\right)  \left(\prod_{j=l+1}^m \overline{v_j}\right) \overline{w}\;dxdt \\&\qquad+\sum_{\substack{I\subset\{1,2,\cdots,m\}\\
\lvert  I \rvert =l,~I\neq \{1,\cdots,l\}}}\displaystyle \int_{Q}
     l!(m-l)!r_l\left(\prod_{j \in I}v_j\right)\left(\prod_{j \notin I}\overline{v_j}\right) \overline{w} \;dxdt\\&\qquad+ \sum_{\substack{a=1\\a\neq l}}^{m-1}\sum_{\substack{I\subset\{1,2,\cdots,m\}\\
\lvert  I \rvert =a}}\displaystyle \int_{Q}
     a!(m-a)!r_a\left(\prod_{j \in I}v_j\right)\left(\prod_{j \notin I}\overline{v_j}\right) \overline{w} \;dxdt.
\end{aligned}
\end{align}
As discussed earlier, we substitute \eqref{vi cgo} and \eqref{w cgo} into the integral identity \eqref{r term}
to arrive at
\begin{align}
\label{integral identity to recover r_l}
\begin{aligned}
    \text{ Known}&=\int_{Q}    m! r_me^{\mathrm{i}(\sum_{k=1}^m\phi_k-\phi_0)} \left(\prod_{k=1}^m U_0^k\right)\overline{U_0^0}~ dx dt\\&\quad+\int_{Q}
     l!(m-l)!r_l e^{\mathrm{i}(\sum_{k=1}^l \phi_k-\sum_{k=l+1}^m\phi_i-\phi_0)}\left(\prod_{k=1}^l U_0^k\right) \left(\prod_{k=l+1}^m \overline{U_0^k}\right)  \overline{U_0^0}\;dxdt \\&\quad+\sum_{\substack{I\subset\{1,2,\cdots,m\}\\
\lvert  I \rvert =l,~I\neq \{1,\cdots,l\}}}\displaystyle \int_{Q}
     l!(m-l)!r_le^{\mathrm{i}\left(\sum_{j\in I}\phi_j- \sum_{j\notin I}\phi_j-\phi_0\right)}\left(\prod_{j \in I}U_0^j\right)\left(\prod_{j \notin I}\overline{U_0^j}\right) \overline{U_0^0} \;dxdt\\&\quad+ \sum_{\substack{a=1\\a\neq l}}^{m-1}\sum_{\substack{J\subset\{1,2,\cdots,m\}\\
\lvert  J \rvert =a}}\displaystyle \int_{Q}
     a!(m-a)!r_ae^{\mathrm{i}\left(\sum_{j\in J}\phi_j- \sum_{j\notin J}\phi_j-\phi_0\right)}\left(\prod_{j \in J}U_0^j\right)\left(\prod_{j \notin J}\overline{U_0^j}\right) \overline{U_0^0} \;dxdt\\&\quad+\mathcal{O}(\lambda^{-1}).
     \end{aligned}
\end{align}
Now, we select $\phi_l$ for $0\leq l\leq m$ so that the exponential factor in the second term of \eqref{integral identity to recover r_l} vanishes while   in other terms of \eqref{integral identity to recover r_l}, it must remain non-vanishing. To achieve this,     
we choose  
$\o_k~(0\leq k\leq m)$ so that 
\begin{align}\label{omega i relation 1}
\begin{aligned}
\o_{1}+\cdots+\o_{l}&=\o_0+\o_{l+1}+\cdots+\o_{{m}},\\ \lvert \o_1\rvert^2+\cdots +\lvert \o_{l}\rvert^2&=\lvert \o_0\rvert^2+\lvert \o_{l+1}\rvert^2+\cdots +\lvert \o_{m}\rvert^2 
\end{aligned}
\end{align}where $1\leq l \leq m-1$. We provide the existence of above mentioned $\omega_k's$ in two different cases viz. $l=1$ and $2\leq l\leq m-1$. \\
\textbf{Case(1)}:
When $l=1$
\begin{table}[ht]
\begin{center}
\begin{tabular}{|Sc|Sc|}
   \hline $m= 2$ & $m\geq 3 $  \\\hline
     $  \o_0=(1,0,\cdots,0)$& $\o_2=\o_3=\cdots=\o_{m}=(0,-1,0,\cdots,0) $\\\hline$\o_2=(0,-1,0,\cdots,0)$&$ \o_0=\left(0,\dfrac{m-2}{2},0,\cdots,0\right)$ \\\hline $\o_1=(1,-1,0\cdots,0)$ &$\o_1=\left(0,-\dfrac{m}{2},0\cdots,0\right)$\\\hline
\end{tabular}
\end{center}
\vspace{0.5cm}
\caption{A particular choice of $\o_k$ for $0\leq l\leq m$ to reconstruct $r_1$}
\label{table2}
\end{table}

\noindent In this case, equation \eqref{omega i relation 1} can be rewritten as

\begin{align}\label{omega i relation l=1}
\begin{aligned}
\o_{1}&=\o_0+\o_{2}+\o_3+\cdots+\o_{{m}},\\ \lvert \o_1\rvert^2&=\lvert \o_0\rvert^2+\lvert \o_{2}\rvert^2+\lvert \o_3\rvert^2+\cdots +\lvert \o_{m}\rvert^2 .
\end{aligned}
\end{align}
A particular choice of $\omega_k's$ satisfying \eqref{omega i relation l=1} can be seen in 
 \Cref{table2}  and with this choice, it is easy  to see that 
\begin{align}
   & \phi_1-(\phi_2+\cdots+\phi_m+\phi_0)\\&= \lambda x\cdot(\o_1-(\o_2+\cdots+\o_{m}-\o_0))- \lambda^2 t(\lvert \o_1\rvert^2-(\lvert \o_2\rvert^2+\cdots \lvert \o_{m}\rvert^2+\lvert \o_0\rvert^2)=0.
\end{align}
Furthermore, for these choices of $\o_k$, the exponential factors other than the second term in \eqref{integral identity to recover r_l} will remain non-vanishing. To prove this, let us consider the following cases:
\begin{enumerate}
    \item If $\sum_{k=1}^m\phi_k-\phi_0=0$ then 
    \begin{align}\label{omega m sum}
         \o_1+\cdots+\o_{m}=\o_0\ \text{ and }\lvert \o_1\rvert^2+\cdots \lvert \o_{m}\rvert^2=\lvert \o_0\rvert^2.
    \end{align}
  The difference of  \eqref{omega i relation l=1} and \eqref{omega m sum} results into 
  \begin{align}
     \o_{2}+\o_3+\cdots+\o_{{m}}=0\ \text{ and }\ \lvert \o_{2}\rvert^2+\lvert \o_3\rvert^2+\cdots +\lvert \o_{m}\rvert^2=0,  
  \end{align}
  which is a contradiction, as $\o_k \in \mathbb{R}^n\setminus\{0\}$ for each $0\leq k\leq m$.
  \item If $\sum_{j\in I}\phi_j- \sum_{j\notin I}\phi_j-\phi_0=0$ for some $I\subset\{1,2,\cdots,m\}$, $\lvert I\rvert=1$ and $I\neq\{1\}$.\\
  Without loss of generality, let us assume $I=\{a\}$, where $a\neq 1$ then we have 
  \begin{align}\label{omega a}
  \begin{aligned}
      \o_a&=\o_1+\o_2+\cdots+\o_{a-1}+\o_{a+1}+\cdots+\o_m+\o_0,\\\lvert    \o_a\rvert^2&=\lvert    \o_1\rvert^2+\lvert    \o_2\rvert^2+\cdots+\lvert    \o_{a-1}\rvert^2+\lvert    \o_{a+1}\rvert^2+\cdots+\lvert    \o_m\rvert^2+\lvert    \o_0\rvert^2.
      \end{aligned}
  \end{align}
  After subtracting \eqref{omega i relation l=1}
 and \eqref{omega a}, we get
 \begin{align}
     \o_a=\o_1 \ \text{ and }\ \lvert  \o_a\rvert^2=\lvert  \o_1\rvert^2.
 \end{align}
 The above expression, together with the second term in \eqref{omega a} leads to
 \begin{align}
     \o_k=0 \text{ for $0\leq k\leq m,~k\neq\{1,a\}$},
 \end{align}
 which is a contradiction. 
 \item If $\sum_{j\in I}\phi_j- \sum_{j\notin I}\phi_j-\phi_0=0$ for some $I\subset \{1,2,\cdots,m\}, \lvert I\rvert=a$ and $a\neq 1$.\\
 Without loss of generality, let $I=\{i_1,\cdots,i_a\}$, then we have
 \begin{align}\label{omega |l|=a}
\begin{aligned}
   \omega_{i_1}+\cdots+  \omega_{i_a}&=\o_{j_1}+\cdots+\o_{j_{m-a}}+\o_0,\\
   \lvert  \omega_{i_1}\rvert^2+\cdots+\lvert  \omega_{i_a}\rvert^2&=\lvert  \omega_{j_1}\rvert^2+\cdots+\lvert  \omega_{j_{m-a}}\rvert^2+\lvert  \omega_{0}\rvert^2.
   \end{aligned}
\end{align}
The difference of \eqref{omega i relation l=1} and \eqref{omega |l|=a} leads to
\begin{enumerate}
    \item Either 
    \begin{align}
        \o_{i_1}+\cdots+\o_{i_{s-1}}+\o_{i_{s+1}}+\cdots+\o_{i_a}=0,\\
        \lvert \o_{i_1}\rvert^2+\cdots+\lvert \o_{i_{s-1}}\rvert^2+\lvert \o_{i_{s+1}}\rvert^2+\cdots+\lvert \o_{i_a}\rvert^2=0,
    \end{align}
    if $i_s=1$ for some $i_s\in I$, which is a contradiction.
    \item Or
    \begin{align}
        \o_1= \omega_{i_1}+\cdots+  \omega_{i_a}\ \text{ and }\ \lvert \o_1\rvert^2=\lvert\omega_{i_1}\rvert^2+\cdots+\lvert  \omega_{i_a}\rvert^2,
    \end{align}
    if $i_s\neq 1$ for any $i_s \in I$. On comparing the above expression  with \eqref{omega i relation l=1}, we get $\o_{j_1}=\cdots=\o_{j_{m-a}}=0$, which is a contradiction.
\end{enumerate}
 
\end{enumerate}
Thus, for these choices of $\o_k$, the exponential factors other than the second term in \eqref{integral identity to recover r_l} will remain non-vanishing. Next,
following a  similar analysis as used in the establishing for reconstruction of $r_m$ in Step \ref{step 4}, we arrive at  
\begin{align}\label{r_1 integral identity}
\begin{aligned}
\text{Known}&=\int_{B(x_0,h)}r_1(t,x)\left(\prod_{k=1}^m \left[\prod_{j=1}^{n-1} \chi\left(
    \dfrac{(x-x_0)\cdot \xi_j^k}{h}
    \right)\right] \prod_{j=1}^{n-1} \chi\left(
    \dfrac{(x-x_0)\cdot \xi_j^0}{h}
    \right)\right)
  \times
  \\&\qquad \bigg[  \exp\bigg(\mathrm{i}\int_{0}^{\infty} \o_1\cdot A(t,x+r \o_1) ~dr -\mathrm{i}\sum_{k=2}^m  \int_{0}^{\infty} \o_k\cdot A(t,x+r \o_k) ~dr\\&\qquad \qquad
  -\mathrm{i} \int_{0}^{\infty} \o_0\cdot A(t,x+r \o_0) ~dr\bigg)\bigg]~dx
    \end{aligned}
\end{align}
for all $t\in (0,T)$. 
Again, we multiply the above expression by $h^{-n}$ and letting $h \rightarrow 0$, to obtain that  
\begin{align}
    r_1(t,x_0)= \text{Known, for all $t\in(0,T)$}.
\end{align}
Now since $x_0 \in \O$ is arbitrary therefore $r_1(t,x)=\text{Known}$, for all $(t,x) \in Q$.\\

\noindent\textbf{Case(2)}: When $2\leq l \leq m-1$\\
Let $\o_k$ for $0\leq k\leq m$ are as in \eqref{omega i relation 1}. Clearly
\begin{align*}
   \sum_{k=1}^l\phi_k - \sum_{k=l+1}^m \phi_k-\phi_0=0.
\end{align*}
Also, except for $J=\{1,\cdots,l\}$, we seek for
\begin{align}\label{sum we want}
\begin{aligned}
    \sum_{j\in J}\phi_j&- \sum_{j\notin J}\phi_j-\phi_0\\&=\lambda x\cdot \left(\sum_{j\in J}\o_j-\sum_{j\notin J}\o_j-\o_0\right)-\lambda^2 t \left(
    \sum_{j\in J}\lvert \o_j\rvert^2-\sum_{j\notin J}\lvert \o_j\rvert^2-\lvert \o_0\rvert^2
    \right)\neq 0
    \end{aligned}
\end{align}
However, for the choice of $\o_k$ for $0\leq k\leq m$ satisfying \eqref{omega i relation 1}, the above expression is not true in general.
\begin{example}
\begin{itemize}
    \item[(1)]  Let $m=3$ and $l=2$, then from \eqref{omega i relation 1}, we have $\o_1+\o_2=\o_3+\o_0$. Choose $\o_1=\o_3=(1,0,\cdots,0)$ and $\o_0=\o_2=(0,1,0,\cdots,0)$. With these choices of $\o_i$, we have \[\phi_2+\phi_3-\phi_1-\phi_0=0.\]
    \item[(2)]Let $m=4$ and $l=2$, then from \eqref{omega i relation 1}, we have $\o_1+\o_2=\o_3+\o_4+\o_0$. Choose $\o_1=\o_3=(1,0,\cdots,0)$, $\o_0=\o_4=(0,1,0\cdots,0)$ and $\o_2=(1,1,0,\cdots,0)$. With these choices of $\o_i$, we have \[\phi_2+\phi_3-\phi_1-\phi_4-\phi_0=0.\]
\end{itemize}
\end{example}
\noindent Now, to ensure that both \eqref{omega i relation 1} and \eqref{sum we want} hold, we have to impose some extra condition on $\o_k$. If possible, let us assume that $\sum_{j\in J}\phi_j- \sum_{j\notin J}\phi_j-\phi_0=0$ for $(J\neq\{1,\cdots,l\}) \subset \{1,2,\cdots,m\}$, then
\begin{align}\label{omega i relation 2}
   \sum_{j\in J}\o_j-\sum_{j\notin J}\o_j-\o_0=0 \quad \text{ and }\quad \sum_{j\in J}\lvert \o_j\rvert^2-\sum_{j\notin J}\lvert \o_j\rvert^2-\lvert \o_0\rvert^2=0.
\end{align}
Subtract \eqref{omega i relation 1} and \eqref{omega i relation 2}, we get
\begin{align}\label{where the contradiction needed}
\sum_{j \in K_1} \o_j=\sum_{j \in K_2} \o_j
\end{align}
where $K_1\subsetneq \{1,\cdots,l\}$ and $K_2 \subsetneq \{l+1,\cdots,m\}$. Next, we 
 choose $\o_k$ for $0\leq k\leq m$ which satisfy \eqref{omega i relation 1} in such a way that \eqref{where the contradiction needed} does not holds for any proper subset $K_1$ and $K_2$, introduced earlier. In particular, we choose $\o_k$ for $0\leq k\leq m$ as in \Cref{table3}.
\begin{table}[ht]
\centering
    \begin{tabular}{|Sc| }
    \hline
 $\o_1=\cdots=\omega_{l-1}=\left(\dfrac{(m-l)(1+m-l)}{l(l-1)}\right)^{\frac{1}{2}}(1,0,\cdots,0)$\\ \hline $\o_l=-(l-1)\left(\dfrac{(m-l)(1+m-l)}{l(l-1)}\right)^{\frac{1}{2}}(1,0,\cdots,0)$ \\\hline
    $\o_{l+1}=\cdots=\o_m=(0,1,0,\cdots,0)$\\\hline $\o_0=-(m-l)(0,1,0,\cdots,0)$\\\hline
    \end{tabular}
    \vspace{0.5cm}
\caption{A particular choice of $\o_k$  for $0\leq k\leq m$ to reconstruct $r_l\ (2\leq l\leq m-1)$.}
\label{table3}
\end{table}

\noindent Again, the foregoing calculations follows from Step \ref{step 4}, and we achieve

\begin{align}\label{r_2 to m integral identity}
\begin{aligned}
\text{Known}&=\int_{B(x_0,h)}r_l(t,x)\left(\prod_{k=1}^m \left[\prod_{j=1}^{n-1} \chi\left(
    \dfrac{(x-x_0)\cdot \xi_j^k}{h}
    \right)\right] \prod_{j=1}^{n-1} \chi\left(
    \dfrac{(x-x_0)\cdot \xi_j^0}{h}
    \right)\right)
  \times
  \\&\qquad \bigg[  \exp\bigg(\sum_{k=1}^l\mathrm{i}\int_{0}^{\infty} \o_k\cdot A(t,x+r \o_k) ~dr -\mathrm{i}\sum_{k=l+1}^m  \int_{0}^{\infty} \o_k\cdot A(t,x+r \o_k) ~dr\\&\qquad \qquad
  -\mathrm{i} \int_{0}^{\infty} \o_0\cdot A(t,x+r \o_0) ~dr\bigg)\bigg]~dxdt
    \end{aligned}
\end{align}
for all $t\in (0,T)$. 
Thus, on multiplying the above expression by $h^{-n}$ and letting $h \rightarrow 0$, we get
\begin{align}
    r_l(t,x_0)= \text{Known,  for all $t\in(0,T)$}.
\end{align}
Since $x_0 \in \O$ is arbitrary therefore $r_l(t,x)=\text{Known}$, for all $(t,x) \in Q$. As $2\leq l \leq m-1$ is arbitrary, we have $r_l(t,x)=\text{Known}$, in $Q$ for every $2\leq l \leq m-1$. This completes the proof.\qed

\section*{Acknowledgments}\label{sec:acknowledgements}
\begin{itemize}
	\item Parveen Kumar acknowledges financial support from the Council of Scientific and Industrial Research (CSIR), India, through the fellowship 09/1005(19269)/2024-EMR I.
    \item Manmohan Vashisth’s work was supported by the ARG-MATRICS grant from the ANRF, Government of India (File No. ANRF/ARGM/2025/002368/MTR).
    \item This research also received partial support under the FIST program of the Department of Science and Technology, Government of India (Ref. No. SR/FST/MS-I/2018/22(C)).\\
\end{itemize}

 \noindent \textbf{Data availability statement.} \ No datasets were generated or analyzed during the current study; therefore, data sharing is not applicable.

    \vspace{.1cm}
\noindent \textbf{Conflict of interest.} \
The authors declare that they have no conflicts of interest regarding the research, authorship, and/or publication of this manuscript.

\end{document}